\documentclass[francais]{svmult}
\usepackage[dvips]{graphicx}
\usepackage[latin1]{inputenc}
\usepackage{multicol}
\usepackage{makeidx}
\usepackage[francais,english]{babel}

\usepackage{amsmath,amsfonts,amssymb}
\usepackage{newlfont}
\usepackage{showidx}

\newtheorem{theo}{Théorème}
\newtheorem{prop}{Proposition}
\newtheorem{cor}{Corollaire}

\newtheorem{defi}{Définition}

\def\g{\mathfrak{g}}\def\a{\mathfrak{a}}
\def\G{\Gamma}
\def\ov{\overline}
\def\tr{\mathrm{tr}}
\def\ad{\mathrm{ad}}
\def\p{\mathfrak{p}}
\def\b{\mathfrak{b}}
\def\h{\mathfrak{h}}

\def\q{\mathfrak{q}}
\def\n{\mathfrak{n}}
\def\k{\mathfrak{k}}

\def\t{\mathfrak{t}}
\def\Exp{\mathrm{Exp}}
\def\d{\mathrm{d}}

\let\ov\overline
\newcommand{\semi}{\,\triangleright\!\!\! <}
\newcommand{\fin}{\hfill $\blacksquare$ \\}
\begin{document}

\bibliographystyle{alpha}
\title*{Applications de la bi-quantification à la théorie de Lie}
\author{ \textbf{Charles Torossian}}
\institute{Université Paris 7, CNRS-Institut Mathématiques
de Jussieu, Équipe de Théorie des Groupes, Case
7012, 2 place Jussieu,  75251 Paris Cedex 05 FRANCE, torossian@math.jussieu.fr}

\maketitle


{\noindent \bf{Abstract :} } This article in French, with a large English introduction, is a survey about applications of bi-quantization theory in Lie theory. We focus on a conjecture of M. Duflo. Most of the applications are coming from our article with Alberto Cattaneo \cite{CT} and some extensions are relating  discussions with my student \cite{BP}. The end of the article is completely new. We prove that the conjecture $E=1$ implies the Kashiwara-Vergne conjecture. Our deformation is non geometric but uses a polynomial deformation of the coefficients.\\

\noindent \textbf{Thanks :} I thank the organizers of the conference ``Higher Structures in Geometry and Physics" which held at IHP, January 2007 and where our results have been  announced. This article is dedicated to  Murray Gerstenhaber and Jim Stasheff.

\section{English Introduction}

\subsection{Invariant differential operators on line bundle}
Let $G$ be a real Lie group, connected and simply connected. Let $\g$ the associated Lie algebra, $U(\g)$ the universal enveloping algebra and  $S(\g)$ the symmetric algebra. In this introduction
 $G/H$ is a homogeneous space with
 $H$ a connected Lie sub-group. As usual $\h$ denotes the Lie algebra of $H$. Fix a character $\lambda$ of $H$, it's a group homomorphism  from $H$ into $\mathbb{C}^{\times}$. If there are no danger of confusion, we will denote by the same letter  the differential of the character. So $\lambda$ is a character of  $\h$, ie. we have $\lambda[\h, \h]=0$.\\

Let's ${\mathcal L}_{\lambda}$ be the line bundle defined by  $\lambda$. Sections of this bundle, denoted by $\Gamma({\mathcal L}_\lambda)$,   are smooth functions on
 $G$ such that
 $\varphi(gh)=\varphi(g)\lambda(h)$. Obviously $G$ acts on the left on $\Gamma({\mathcal L}_\lambda)$.\\

Let ${\mathcal D}_{\lambda }$ be the algebra of invariant differential operators on $\Gamma({\mathcal L}_\lambda)$.\\

After Koornwinder \cite{koorn}  we know ${\mathcal D}_{\lambda }$  is isomorphic to
\begin{equation}{\mathcal D}_{\lambda }:= \Big(U(\g)_\mathbb{C} /U(\g)_\mathbb{C}\cdot \h_{-\lambda}\Big)^{\h}\end{equation}
 where
$\h_{-\lambda}=\{H-\lambda(H) \}$. Here are some explanations. For $X\in \g$, $R_X$ is the left invariant vectors field on $G$ associated to $X$. For
$u \in {\mathcal D}_{\lambda }$, let  $D_{u}$ be
 the associated differential  operator on $\Gamma({\mathcal L}_\lambda )$
defined by
$$(D_{u}\varphi)(g)=(R_{u}\varphi)(g),$$$\varphi \in \Gamma({\mathcal L}_\lambda)$. Then $D_u$ in a left invariant differential operator. It's not difficult to verify that we have described all of them\footnote{Consider the local expression around the origin.}.

Suppose $\lambda$ real. The algebra $\Big(U(\g)/U(\g)\cdot\h_{-\lambda}\Big)^\h$  is not commutative in general. A conjecture of M. Duflo \cite{duflo-japon} describes the center of this algebra. I write $\mathcal{S}_\lambda$ the algebra of $H$-invariant polynomial functions on $ \h_{-\lambda}^\perp :=\{ f\in\g^*, f|_\h=\lambda \}$.  We get  \begin{equation}\mathcal{S}_\lambda :=
  \Big(S(\g)/S(\g)\cdot \h_{-\lambda}\Big)^\h.\end{equation}
  This space admits a natural Poisson structure coming from the classical Poisson structure on $\g^*$. Put $\delta(H)=\frac12 \tr_{\g/\h} \, \ad (H)$ the character for the half densities.\\

\noindent \textbf{Duflo's conjecture \cite{duflo-japon} : }The center of  $\Big(U(\g)/U(\g)\cdot\h_{-\lambda-\delta}\Big)^\h$ is isomorphic to the Poisson center of
  $\mathcal{S}_\lambda$.\\

This conjecture is far to be solved\footnote{Consider the following example ;  $G$
is reductive $H=U$ the unipotent radical of a Borel. Let's note  $\t$ a Cartan subalgebra. Then you get $(U(\g)/U(\g)\cdot
\h)^{\h}=S(\t)=(S(\g)/S(\g)\cdot \h)^{\h}$. These algebras are commutative. The space   $G/U$ is quasi-affine.}. Moreover one should probably ask for generic character.\\

 In case
$G$ is a  nilpotent group, appreciable advances have been achieved in last few years by Corwin-Greenleaf \cite{CG},
Fujiwara-Lion-Magneron-Mehdi \cite{FLMM}, Baklouti-Fujiwara
\cite{BF} and  Baklouti-Ludwig \cite{BL} and Lipsman \cite{Li1, Li2, Li3}. More precisely in the nilpotent case one can prove the following.

\begin{theo}[\cite{FLMM}] Let  $G$ nilpotent (connected, simply connected) and  $\chi$ the unitary character of $H$ defined by   $\chi(\exp_G(H))=\exp(i\lambda(H))$.
\begin{multline}\nonumber \mathcal{D}_\chi \; \mathrm{commutative} \Longleftrightarrow \mathrm{the \; left \;  representation \; } L^2(G/H, \chi) \; \mathrm{has\;finite \; multiplicities }\\
\Leftrightarrow \mathrm{for\; } f \in \h_{-\lambda}^\perp \; \mathrm{generic \;}  \quad
H\cdot f  \mathrm{ \; is\; \; lagrangian  \;in  \; }  G\cdot f\\\Leftrightarrow
\Big( \mathrm{Frac} (S(\g)/S(\g)\cdot \h_{-\lambda}) \Big)^\h \quad
\mathrm{\;Poisson\; commutative}\\ \Leftrightarrow \Big(S(\g)/S(\g)\cdot
\h_{-\lambda}\Big)^\h \quad \mathrm{ Poisson \; commutative}
\end{multline}
Under these conditions  $\mathcal{D}_\chi$ is a subalgebra of  $\Big(
\mathrm{Frac} (S(\g)/S(\g)\cdot \h_{-\lambda}) \Big)^\h$, the algebra of $H$-invariant fractions on $\h_{-\lambda}^\perp$.
\end{theo}

In case where $G$ and $H$ are reductive groups  F. Knop \cite{knop} gives a satisfying  and remarkable answer to the conjecture. In case $H$ is compact and  $G=H\semi N$ is a semiproduct of $H$ with a
Heisenberg group $N$,  L. Rybnikov
 \cite{ryb} makes use of Knop's result to prove the Duflo's conjecture. The case of  symmetric spaces has been previously studied by myself in  \cite{To3, To4} and the group case is solved by Duflo \cite{du77}.\\\\

This survey analyses  the  Duflo conjecture and other standard problems in Lie theory with the help of Kontsevich's quantization.  Main results are some generalizations of  \cite{CT}. Some others comes from discussion with my PhD student P. Batakidis \cite{BP}, the end of the article  is new.\\

One hopes to convince people of the interest of our methods. In some sense they are a replacement for the orbit method.\\

\subsection{Duflo's conjecture  : a review of difficulties }
Let's try to list some technical difficulties in Duflo's conjecture. We will see that most of them disappear with Kontsevich's quantization techniques developed   in  \cite{CT}.\\

1 - The algebra $\Big(U(\g)/U(\g)\cdot \h_{-\lambda}\Big)^{\h}$ is filtered by the order of differential operators. But it's not obvious to describe the associated graded space. In general we have a  injection  from
$gr\Big(U(\g)/U(\g)\cdot \h_{-\lambda}\Big)^{\h}$ into
 $\Big(S(\g)/S(\g)\cdot \h\Big)^{\h}$. Let's remark that the character has  disappeared.  It's the first difficulty: the symbol of a differential operator on a line bundle is just a function on the cotangent space of the underlying space.  Except tentative  \cite{Kostant} there are no way to keep the character. The next example illustrates the phenomena : consider   $\g=< X, Y, Z >$ with $Z=[X, Y]$ and
    $\h=< Z >$.\\

 If  $\lambda(Z)=0$ the algebra
    $(U(\g)/U(\g)\cdot Z)^Z=S(\g)/S(\g)\cdot Z$ is Poisson commutative.

    If  $\lambda(Z)\neq 0$ then  $\Big(U(\g)/U(\g)\cdot
    (Z-\lambda(Z))\Big)^Z$ is a Weyl algebra then non commutative while the associated graded is  Poisson
    commutative.  There are no link between  ${\mathcal D}_{\lambda }$ and the  Poisson algebra
    $gr({\mathcal D}_{\lambda })$. The quantization procedure developed by Cattaneo-Felder  \cite{CF1, CF2} for  co-isotropic space takes care of the character.\\

2 - In general the homogeneous space $G/H$ doesn't admit a  $G$-invariant measure. One has to consider half-densities, essentially to deal with Hilbert spaces. So you have to define the following character of $H$, $\Delta_{G,H}(h)=(\det_{\g/\h} \mathrm{Ad} h)^{1/2}$ and  $\lambda$ is replaced by the shifted character
$\lambda +\frac12 \tr_{\g/\h} \ad$. This shift is problematic if you want to construct irreducible representations by induction from a polarization $\b$. Usually you have to ask for  compatibility conditions among $\Delta_{G,H}, \Delta_{G,B}$ \cite{Li1, To3}. Our theories  don't use  this shift; they are well normalized.\\

3 - In case  $\h$ admits a $\h$-invariant complement (for the adjoint action), the pair  $(\g, \h)$ is called a \textit{reductive pair} (but $\g$ is not supposed to be reductive !). If $(\g, \h)$  is a  reductive pair, then an easy consequence of the   Poincaré-Birkhoff-Witt's theorem indicates $\big(U(\g)/U(\g)\cdot \h_{-\lambda} \big)^\h$ is isomorphic as a  vector space to
$\big(S(\g)/S(\g)\cdot \h_{-\lambda}\big)^\h$. It is not know whether  this holds in  general, that means if  $\mathcal{D}_\lambda$ is a deformation of  $\mathcal{S}_\lambda$.  Actually there are no obvious map from
 $\big(S(\g)/S(\g)\cdot \h_{-\lambda}\big)^\h$ into $U(\g)$ or  $U(\g)/U(\g)\cdot \h_{-\lambda}$.\\

 If  $(\g, \h)$ is a reductive pair then you get the extra equality
\begin{equation}U(\g)^{\h}/U(\g)^{\h}\cap U(\g)\cdot \h_{-\lambda}=
(U(\g)/U(\g)\cdot \h_{-\lambda})^{\h}.\end{equation}
In general you get just a injection from LHS to RHS, as illustrates the following example. Consider  $\g=\textbf{sl}(2)$ with standard basis $H, X, Y$ and take
 $\h= < X >$. You get
$\big(U(\g)/U(\g)\cdot X\big)^{X}=\mathbb{R}[X]$ but  $U(\g)^{X}/U(\g)^{X}\cap U(\g)\cdot X$ is isomorphic to $\mathbb{R}[X^2]$. Our constructions depend of the choice for a complement to $\h$.  In \cite{CT} we gave several examples where there are different choices for a complement, the most important is the Iwasawa decomposition and the Cartan decomposition for symmetric pairs.  Dependance of our constructions with the complement leads to interesting applications : Harish-Chandra homomorphism for example. At the end of this article, we will explain the group case\footnote{This is the double $(\g\times \g, \mathrm{diagonal})$.}, which leads to the Kashiwara-Vergne conjecture. Invariant complements simplify the calculations but we are able to describe our model even in the general case.\\

4 - The algebras  $\big(U(\g)/U(\g)\cdot \h_{-\lambda} \big)^\h$ and $\mathcal{S}_\lambda$ should be simultaneously commutative. Even this fundamental question  is not solved in general, except for nilpotent case or for symmetric pairs. If  $G$ and $H$  are algebraic and the generic  $H$-orbits in
$\h_{-\lambda}^{\perp}=\{ f \in \g^{*}, f|_{\h}=\lambda\}$ are lagrangian then $\Big(
\mathrm{Frac} (S(\g)/S(\g)\cdot \h_{-\lambda}) \Big)^\h$ is commutative,  consequently
 $\mathcal{S}_\lambda$ is Poisson commutative. Of course $\mathcal{S}_\lambda$  could be commutative, without the lagrangian hypothesis. For example consider $\g=< T, X, Y, Z>$ with $< X, Y, Z >$
a Heisenberg Lie algebra and   $[T, X]=X, [T, Y]=Y $  and  $[T, Z]=2Z$. Take $\h=< T >$.
Then  $\big(U(\g)/U(\g)\cdot \h\big)^\h = \mathbb{R}$, $\big(S(\g)/S(\g)\cdot \h\big)^\h=\mathbb{R}$ but $H$-orbits  in $\h^\perp$ are not lagrangian because  $\mathrm{Frac}(S(\g)/S(\g)\cdot \h)^\h$ is not commutative.\\

5 - Fix $\q$ a complement of $\h$ in $\g$ and consider the Exponential map  $\Exp : \q \longrightarrow G/H$
defined by $X \mapsto \exp_G(X) H$. This is a local diffeomorphism and you can write $D_u$ for  $u \in \mathcal{D}_\lambda$ in exponential coordinates. Before our work  \cite{CT} no formulas were known. If you restrict these operators to invariant distributions, you should get interesting simplifications. This is exactly what happens for symmetric pairs
\cite{Rou86, Rou90, Rou91}, especially for the double
$G/H=G_1\times G_1/\mathrm{Diagonal}$. The study of this restriction gives rise to Kashiwara-Vergne's conjecture \cite{KV, AM2, AST, To5}.\\

6 - Suppose $\chi=i\lambda$ is the differential of a unitary character of $H$ and  $\mathcal{D}_\chi$  commutative.
 How can we associate to $u \in \mathcal{D}_\chi$ a rational function, or a polynomial function on
 $\h_{-\lambda}^{\perp}$ ? The orbit method gives a kind of answer : construct an irreducible representation
  $(\pi, \mathcal{H})$ of $G$ which admits
    $H$-semi-invariant  distribution vectors for  $\chi$. Most of them are related with orbits $\Omega=G\cdot f$ with   $f\in \h_{-\lambda}^{\perp}$. If we are lucky, these $H$-semi invariant  distribution vectors are commun eigenvectors for all $D_u$. The eigenvalue is a character for $\mathcal{D}_\lambda$ and  should  depend on $f$ as a rational function. Usually an  irreducible representation $(\pi, \mathcal{H})$ is constructed by induction from a polarization at $f$ (if such a polarization exists !). As we see, there are several analytic difficulties : definition of the  distribution vector, $L^2$ convergence, real structure. All these problems are in some sense far from our starting algebraic problem. We will explain how the bi-quantization gives us a systematic procedure, under the lagrangian hypothesis, to construct this character in a more algebraic (or geometric) way. Of course, if you deal with the spectral decomposition of  $L^2(G/H, \chi)$, all these problems are to be considered. \\

7 - Let's do a fundamental remark now : from the point of view of the theory of representations, one should study  the algebra $(U(\g)/I)^\g$ where $I$ is a  two-sided ideal included in  $U(\g)\cdot \h_{-\lambda}$ and maximal. This algebra should be smaller than  $\Big(U(\g)/U(\g)\cdot\h_{-\lambda}\Big)^\h$ and behaves in a much  better matter. Indeed Schur's lemma proves that the action of any element $u\in (U(\g)/I)^\g$ is scalar on irreducible representations which  admit $H$-semi-invariant  distribution vectors.  It is not difficult to extend Duflo's arguments   \cite{du77} in this context. The rational function you should have built by the orbit method (if it exists!) is then a polynomial function.


\section{La quantification de Kontsevich}
Pour simplifier la compréhension de cet article on rappelle brièvement les constructions de Kontsevich \cite{Kont} et les extensions dans le cas co-isotrope dûes à Cattaneo-Felder \cite{CF1, CF2}
\subsection{Théorème de Formalité}
En 1997, M. Kontsevich a montré que tout variété de Poisson admet une quantification formelle. C'est une conséquence du théorème de formalité qui affirme qu'il existe un quasi-isomorphisme entre l'algèbre de Lie des polychamps de vecteurs munie du crochet de Schouten et l'algèbre de Lie des opérateurs polydifférentiels munie du crochet de Gerstenhaber \cite{Ger} et de la différentielle de Hochschild.

\begin{theo}[\cite{Kont}]\label{theoKont} Il existe un $L_\infty$-quasi isomorphisme $\mathcal{U}=(U_n)_{n\geq 1}$ entre les algèbres différentielles graduées
 $\g_1=(T_{poly}(\mathbb{R}^d, [\cdot , \cdot ]_S, \d=0)$ et $
\g_2=D_{poly}(\mathbb{R}^d, [\cdot, \cdot]_{G}, \d_{Hoch} )$.  En particulier $\mathcal{U}$ induit une bijection entre les solutions formelles de Maurer-Cartan modulo les groupes de jauge.
\end{theo}

La preuve du théorème utilise une construction explicite en terme de diagrammes, pour décrire les coefficients de Taylor $U_n$ de $\mathcal{U}$. En particulier si  $\pi$ est un bi-vecteur  de Poisson vérifiant $[\pi, \pi]_S=0$ pour le crochet de Schouten, alors
\begin{equation}\underset{\epsilon}{\star}= m+
\sum_{n\geq 1} \frac{\epsilon^n}{n!}U_n(\underset{n \,
\mathrm{fois}}{\underbrace{\pi, \ldots, \pi}})\end{equation}  est
une structure associative formelle sur $\mathcal{C}^\infty(\mathbb{R}^d)[[\epsilon]]$. Pour $f, g$ des fonctions de $\mathcal{C}^\infty(\mathbb{R}^d)$ on obtient un produit formel associatif

\begin{equation}   f\underset{Kont}\star g=fg +\sum
_{n=1}^{\infty}\frac{\epsilon^{n}}{n!}\sum_{\Gamma\in G_{n,2} }w_{\Gamma}
B_{\Gamma}(f,g).
\end{equation}

On explique rapidement la signification de chaque termes de cette formule.\\

\subsubsection{Graphes } Ici $G_{n,2}$ désigne l'ensemble des graphes
étiquetés\footnote{Par graphe étiqueté
on entend un graphe $\Gamma$ muni d'un ordre total sur l'ensemble
$E_\Gamma$ de ses ar\^etes, compatible avec l'ordre des sommets.}  et orientés (les arêtes sont orientées) ayant $n$ sommets de première espèce numérotés $1,2,
\cdots,  n$
et deux sommets de deuxième espèce $\overline{1}, \overline{2}$,  tels que~:

\smallskip
i-  Les arêtes partent des sommets de première espèce. De chaque sommet de première espèce partent exactement deux arêtes.

\smallskip
ii-  Le but d'une ar\^ete est différent de sa source (il n'y a pas
de boucle).

\smallskip
iii- Il n'y a pas d'ar\^ete multiple. \\

\begin{figure}[!h]
\begin{center}
\includegraphics[width=7cm]{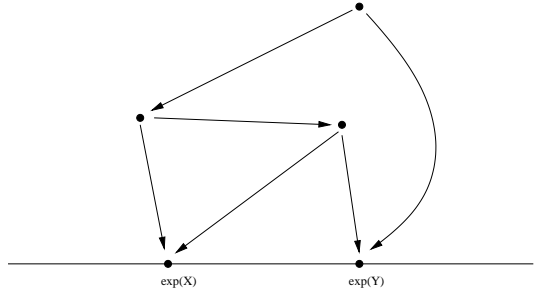}
\caption{\footnotesize Graphe simple de type Lie et de symbole $\G(X, Y)=[[X, [X, Y]], Y]$.}\label{Lie.eps}
\end{center}
\end{figure}

\begin{figure}[h!]
\begin{center}
\includegraphics[width=7cm]{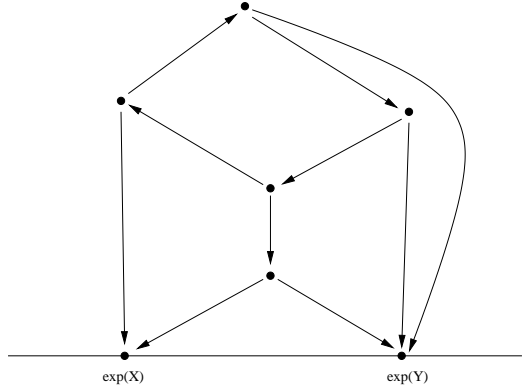}\caption{\footnotesize
  Graphe de type roue  et de symbole $ \Gamma(X,Y)=\tr_\g( \ad X \ad[X,Y]\ad Y\ad
Y).$}\label{roue.eps}
\end{center}
\end{figure}

\noindent \textbf{Remarque importante : } Dans le cas linéaire qui nous intéresse, les graphes qui interviennent de manière non triviale (on dira essentiels), sont tels que les sommets de première espèce ne peuvent recevoir qu'au plus une arête. Il en résulte que tout graphe essentiel est superposition de graphes simples de type Lie (graphe ayant une seule racine comme dans Fig.~\ref{Lie.eps}) ou de type roue (cf. Fig.~\ref{roue.eps} pour un exemple). Cela implique que toutes les formules sont des exponentielles. \\

\subsubsection{Variétés de configurations }On note $C_{n,m}$ l'espace des configurations de $n$ points distincts
dans le
demi-plan de Poincaré (points de première espèce ou points aériens ) et $m$ points  distincts sur la droite réelle
(ce sont les points de seconde espèce ou points terrestres ), modulo l'action
du groupe $az+b$ (pour $a \in \mathbb{R}^{+*}, b\in \mathbb{R}$).  Dans son article \cite{Kont} Kontsevich construit des
compactifications de ces variétés notées $\overline{C}_{n,m}$. Ce
sont des variétés à coins de dimension $2n-2+m$. Ces variétés ne sont pas connexes pour $m\geq 2$.
On notera par  $\overline{C}^{+}_{n,m}$ la composante  qui contient les configurations où les
 points terrestres sont ordonnés dans l'ordre croissant (ie. on  a $\ov{1}< \ov{2}<\cdots< \ov{m}$).

\subsubsection{Fonctions d'angle et coefficients }On définit la fonction d'angle hyperbolique dans le demi-plan de Poincaré par
\begin{equation}
\phi(p,q)=\arg(p-q) + \arg(p-\ov{q}).
\end{equation}
C'est une fonction d'angle de $C_{2,0}$ dans $\mathbb{S}^1$ qui s'étend en une fonction régulière à la compactification  $\overline{C}_{2,0}$. Si $\Gamma$ est un graphe dans $G_{n,2}$, alors  toute arête $e$ définit par restriction
une fonction d'angle notée $\phi_{e}$ sur la variété $\overline{C}^+_{n,2}$. On note $E_{\Gamma}$
l'ensemble des arêtes du graphe $\Gamma$. Le produit ordonné
\begin{equation}
\Omega_{\Gamma}=\bigwedge _{e \in E_{\Gamma}}\d\phi_{e}
\end{equation}
est donc une $2n$-forme sur  $\overline{C}^+_{n,2}$ variété compacte de dimension $2n$.
Le poids associé à un graphe $\G$ est par définition

\begin{equation}
w_{\Gamma}=\frac{1}{(2\pi)^{2n}}\int_{\overline{C}^+_{n,2}} \Omega_{\Gamma}.
\end{equation}
\subsubsection{Opérateurs bi-différentiels }Enfin l'opérateur $B_\Gamma$ est un opérateur bidifférentiel construit à partir de $\Gamma$, dont on ne détaille pas la construction. Disons que chaque arête correspond à une dérivée, chaque sommet de première espèce est attaché au bi-vecteur de Poisson et chaque sommet de deuxième espèce est attaché  à des fonctions (cf. \cite{Kont, CKT}).\\

\subsection{Formule de Baker-Campbell-Hausdorff }

On applique ce théorème pour $\mathbb{R}^d=\g^*$ et
$\pi=\frac12 \sum_{i, j} [e_i, e_j] \partial_{e_i^*}\wedge
\partial_{e_j^*}$
Prenons maintenant $X, Y \in \g$ et $f=e^X, g=e^Y$. L'équation  ci-dessus donne alors une expression nouvelle pour la formule de Baker-Campbell-Hausdorff $Z(X,Y)$; elle utilise tous les crochets possibles \footnote{La formule de Dynkin n'utilise que des crochets itérés.} \cite{Ka}

\begin{equation}Z(X,Y)=X+Y
+\sum\limits_{n\geq 1}\sum\limits_{\substack{
 \Gamma  \;\mathrm{simple}\\\mathrm{geometric}\\
\mathrm{Lie\; type\;  (n,2)}}}w_{\Gamma} \Gamma(X,Y).
\end{equation}

Le terme  $\Gamma(X, Y)$ est le mot de type Lie que l'on peut fabriquer
avec $\Gamma$, c'est essentiellement  le symbole de l'opérateur $B_\G$.

\section{La quantification de  Cattaneo--Felder}

Soit $\g$ une algèbre de Lie de dimension finie sur $\mathbb{R}$.
L'espace dual $\g^*$ est alors muni d'une structure de Poisson
linéaire. On note
 $\pi$ le bi-vecteur de Poisson associé à la moitié du crochet de Lie. Supposons donnés $\h\subset\g$ sous-algèbre de $\g$ et $\lambda$ un caractère réel de $\h$, c'est à dire une forme linéaire telle que $\lambda[\h, \h]=0$.  L'orthogonal
$\h^\perp$ de même que $\h_{-\lambda}^\perp :=\{ f\in\g^*, f|_\h=\lambda \}$ sont des sous-variétés coisotropes de $\g^*$. \\

\subsection{Construction par transformée de Fourier impaire \cite{CF2}}
La construction de \cite{CF1, CF2} concerne le cas des variétés co-isotropes en général, mais nous ne nous intéressons ici qu'au cas des sous-algèbres d'une algèbre de Lie. Les constructions sont locales et dépendent donc d'un choix d'un sup\-plémentaire de $\h$ dans $\g$. Notons $\q$ un tel
supplémentaire. On peut alors identifier $\h^*$ avec $\q^\perp$. On aura une décomposition (affine) $$\g^*= \h_{-\lambda}^\perp\oplus \q^\perp= \h_{-\lambda}^\perp\oplus \h^*.$$
La variété qui intervient dans cette construction est une
super-variété intrinsèque  :

\begin{equation}M :=\h_{-\lambda}^\perp\oplus \mathrm{\Pi} \h
\end{equation}
où $\mathrm{\Pi}$ désigne le
foncteur de  changement de parité. L'algèbre des fonctions polynomiales\footnote{Comme on ne considère que des structures de Poisson linéaires, on peut restreindre les constructions aux algèbres de fonctions polynomiales.} est donc canoniquement
\[\mathcal{A} := \mathrm{Poly}(\h_{-\lambda}^\perp)\otimes
\bigwedge(\g^*/\h^\perp)\simeq \Big(S(\g)/S(\g)\cdot\h_{-\lambda}\Big)\otimes
\bigwedge \h^*.\]
Considérons alors $\pi$ le bi-vecteur de Poisson linéaire de $\g^*$ et appliquons la transformée de Fourier impaire \cite{CF2} dans la direction normale $\h^*=\q^\perp$. On obtient un polyvecteur $\widehat{\pi}$ sur $M$, solution de l'équation de Maurer-Cartan  $$[\widehat{\pi}, \widehat{\pi}]_S=0.$$
On applique  le théorème de formalité pour la super-variété $M$ (voir \cite{CF2} pour une description complète du théorème de Formalité dans le cas gradué),  on obtient alors une solution $\mu$ de Maurer-Cartan dans $\mathcal{D}_{poly}(\mathcal{A})$; c'est un opérateur
polydifférentiel formel homogène de degré un  si l'on
tient compte des degrés impairs.\\

\noindent En d'autres termes, comme la graduation tient compte du degré dans
les variables impaires, la structure obtenue est en fait une
$A_\infty$-structure, définie par Stasheff \cite{Sta}, sur l'espace
$\mathcal{A}=\mathrm{Poly}(\h_{-\lambda}^\perp)\otimes \bigwedge \h^*$
avec premier terme non nul \textit{a priori} (c'est l'anomalie), c'est à dire une
structure
\begin{equation}\mu=\mu_{-1}+ \mu_0+ \mu_1+ \mu_2+ \ldots \end{equation}
vérifiant\footnote{On note $[\bullet,\bullet]_{G}$ le crochet de
Gerstenhaber.} $\frac 12[\mu,
\mu]_{G}=0$ et $\mu_i$ des opérateurs $(i+1)$-polydifférentiels. \\

Dans le cas linéaire, l'anomalie $\mu_{-1}$ est nulle, par conséquent $\mu_0$ est une différentielle et $\mu_1$ un produit associatif modulo des termes contenant $\mu_0 $ et $\mu_2$. Le terme de plus bas degré de $\mu_1$ correspond au crochet de Poisson. Ainsi on construit un vrai produit associatif sur l'espace de cohomologie défini par~$\mu_0$.

\begin{defi} On notera $H^{\bullet}(\mu_0, \mathcal{A})$ l'algèbre de cohomologie (graduée) munie de sa loi associative $\mu_1$. On s'intéressera à la sous-algèbre en  degré $0$ que l'on appellera \textit{algèbre de réduction} et que l'on notera $H^0(\mu_0)$. Le produit $\mu_1$ se restreint en un star-produit  noté $\underset{CF}\star$.
\end{defi}

\subsection{Construction en termes de diagrammes de Feynman \cite{CF1}}\label{sectionrappels}
\begin{figure}[h!]
\begin{center}
\includegraphics[width=8cm]{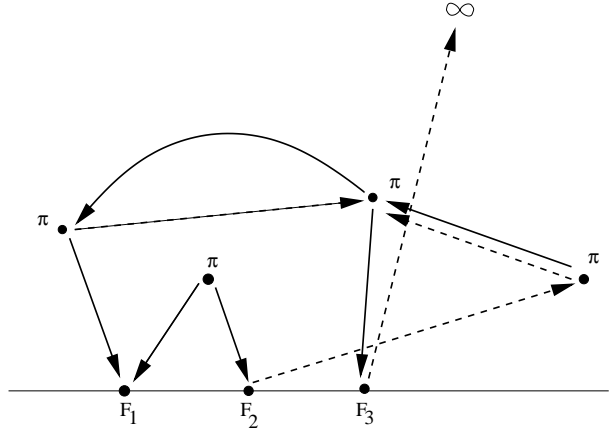}
\caption{\footnotesize Graphe type intervenant dans le calcul de
$\mathcal{U}_4(\pi, \pi, \pi, \pi) $ pour un bivecteur $\pi$ non linéaire. }\label{figuregraphetype}
\end{center}
\end{figure}

La formule proposée est semblable à celle de Kontsevich \cite{Kont}
dans $\mathbb{R}^n$. Chaque $\mu_i$, opérateur $(i+1)$-polydifférentiel s'exprime sous la forme\footnote{Il n'y a pas de terme pour $n=0$ sauf pour $\mu_1$ où on trouve la multiplication $m$.}

\begin{equation}\mu_i= \sum\limits_{n\geq 0} \frac{\epsilon^n}{n!} \sum\limits_{\G \in G_{n, i+1}} w_\G B_\G.\end{equation}
où les arêtes des graphes portent deux couleurs\footnote{Chaque couleur indiquera si la variable de dérivation est dans $\h^*$ ou $\h_{-\lambda}^\perp$ et précisera la fonction d'angle.}. Chaque $B_\G$ est un opérateur $(i+1)$-polydifférentiel sur $\mathcal{A}=\mathrm{Poly}(\h_{-\lambda}^\perp)\otimes \bigwedge \h^*$. Il faut donc élargir la notion de graphes admissibles et considérer des graphes avec arêtes colorées par $\h^*$ issus des points terrestres. Si le bi-vecteur de Poisson $\pi$ n'est pas linéaire, ces graphes peuvent admettre   des arêtes doubles si elles ne portent pas la même couleur. On ne conservera que $2n+ i-1$ arêtes (la dimension de la variété $C_{n, i+1}$) les arêtes restantes seront colorées par $\h^*$ (on dira que ces arêtes vont à l'infini), elles ne contribuent pas dans le calcul du coefficient $w_\G$, mais  les arêtes qui partent à
l'infini contribuent dans la définition de l'opérateur $B_\G$ (cf. Fig. \ref{figuregraphetype}).\\

Concernant le coefficient $w_\G$, il est obtenu de manière similaire par intégration sur la variété $\ov{C_{n, i+1}}^+$ de la forme $\Omega_\G$ modifiée par la couleur selon les règles suivantes :

- si la couleur
est dans $\h_{-\lambda}^\perp$  (variable tangente)  la fonction d'angle
associée est la même que dans le cas classique
\begin{equation}\mathrm{d} \phi_{+}(p,q) :=  \mathrm{d}
\overset{\longrightarrow}{\phi}(p,q):= \mathrm{d}\arg(p-q) +
\mathrm{d} \arg(p-\overline{q}).\end{equation}

- si la couleur est dans $\g^*/\h^\perp=\h^*$ (variable normale\footnote{
On a besoin ici de faire un choix d'un sup\-plémentaire de $\h$,
pour identifier $\h^*$ à un sous-espace de~$\g^*$.})    alors la
fonction d'angle sera notée $\dashrightarrow$ (en pointillé dans les
diagrammes).
\begin{equation} \mathrm{d} \phi_{-}(p,q) := \mathrm{d}
\overset{\dashrightarrow}{\phi}(p,q):= \mathrm{d}\arg(p-q) -
\mathrm{d} \arg(p-\overline{q}).\end{equation}

\subsection{Description de la différentielle $\mu_0$ et exemples d'algèbres de réduction}\label{sectionexemple}
Dans \cite{CT} on décrit en termes de diagrammes ce que vaut la différentielle $\mu_0$ sur les fonctions polynomiales\footnote{On renvoie à \cite{CT} pour l'action sur les éléments de $\mathrm{Poly}(\h_{-\lambda}^\perp)\otimes \bigwedge \h^*.$} sur $\h_{-\lambda}^\perp$.

\begin{figure}[h!]
\begin{center}
\includegraphics[width=6cm]{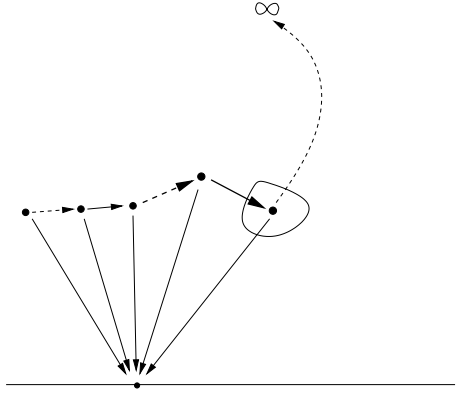}
\caption{\footnotesize Graphe de type Bernoulli}\label{bernoulli2.eps}
\end{center}
\end{figure}

 \begin{figure}[h!]
\begin{center}
\includegraphics[width=6cm]{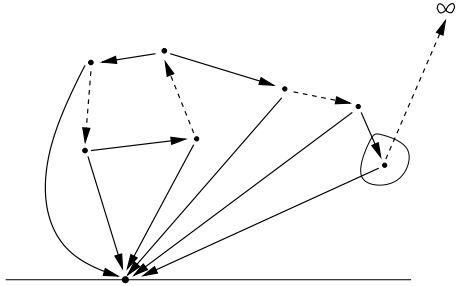}
\caption{\footnotesize Graphe de type roue attaché à un Bernoulli}\label{rouebernoulli.eps}
\end{center}
\end{figure}

\begin{figure}[h!]
\begin{center}
\includegraphics[width=6cm]{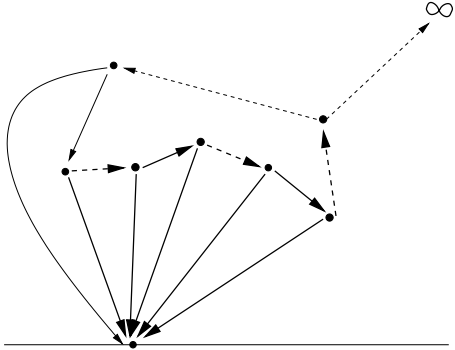}
\caption{\footnotesize Graphe de type roue pure}\label{roueinfini.eps}
\end{center}
\end{figure}
\begin{prop}[\cite{CF1, CT}]  La différentielle $\mu_0$ est l'action de tous les graphes de types suivant :

-i-  les graphes de type Bernoulli avec la dernière ar\^ete partant
\`a
l'infini (cf. Fig.~\ref{bernoulli2.eps})

-ii- les graphes de type roues avec des rayons attachés directement
\`a l'axe réel sauf pour l'un d'entre eux qui est attaché \`a un
graphe de type Bernoulli dont la dernière ar\^ete part \`a l'infini
 (cf. Fig.~\ref{rouebernoulli.eps})

-iii- les  graphes de type roues avec des rayons attachés
directement \`a l'axe réel sauf pour l'un d'entre eux qui part \`a
l'infini
 (cf. Fig.~\ref{roueinfini.eps}).

 En particulier on a toujours $\mu_0= \epsilon \d_{CH} + O(\epsilon^2)$ avec $\d_{CH}$ la différentielle de Cartan-Eilenberg.
 \end{prop}

Donnons quatre exemples d'algèbres de réduction (voir \cite{CT} \S 2 pour les détails).\\

-- Soit $0^\perp=\g^*$. Alors l'algèbre de réduction est $S(\g)$ muni du produit de Kontsevich. En effet il n'y a pas d'arêtes sortantes, donc pas de condition.\\

-- Supposons que $\h$ admette un supplémentaire stable $\q$ alors on montre que $\mu_0= \epsilon \d_{CH}$. Dans ce cas, on en déduit que l'algèbre de réduction s'identifie à $\Big( S(\g)/S(\g)\cdot \h_{-\lambda}\Big)^\h$. En effet, les graphes avec plus de deux sommets ont tous  deux arêtes d'une même couleur qui se suivent (cf. Fig.  \ref{bernoulli2.eps}). Le coefficient associé est alors $0$.\\

-- Soit $f\in \g^*$ et soit $\b$ une polarisation en $f$, c'est à dire une sous-algèbre subordonnée $f[\b, \b]=0$ et lagrangienne pour $B_f(x,y)=f[x,y]$. On prend  comme espace affine $f + \b^\perp$, alors l'algèbre de réduction vaut $\mathbb{C}$.\\

-- Soit $\g=\k  \oplus \a \oplus \n$ une décomposition d'Iwasawa d'une algèbre réductive réelle. On note $\mathfrak{m}$ le centralisateur de $\a$ dans $\k$ et on prend comme sous-algèbre $\mathfrak{m} \oplus \n$, alors l'algèbre de réduction vaut $S(\a)$.

\paragraph{Parité :} On montre \cite{CT} \S 2.2 pour $\lambda =0$, que dans la différentielle $\mu_0$ seuls interviennent les diagrammes avec un nombre impair de sommets de première espèce. En tenant compte du degré de $\epsilon$, on en déduit que l'algèbre de réduction est graduée. Toute fonction $F$ homogène de degré total $n$ dans l'algèbre de réduction  s'écrit $$F= F_n + \epsilon^2 F_{n-2} + \epsilon^4 F_{n-4} + \ldots,$$
avec $F_i$ de degré $i$. On a $F_n \in \Big( S(\g)/S(\g)\cdot \h\Big)^\h;$ c'est l'analogue du symbole. Toutefois on ne peut en déduire que l'algèbre de réduction est une quantification de $\Big( S(\g)/S(\g)\cdot \h\Big)^\h$ car l'application \textit{limite classique} $LC : F\mapsto F_n$ est injective mais pas forcément surjective. Dans le cas général $\lambda \neq 0$ on peut étendre la construction ci-dessus \cite{BP}. On déduit classiquement le corollaire suivant.

\begin{cor}\label{com} Si l'algèbre de réduction $H^0(\mu_0)$  est commutative, alors  son image par l'application LC  est Poisson commutative.
\end{cor}

\subsection{Bi-quantification de Cattaneo-Felder}
On suppose données  deux sous-algèbres $\h_1$ et $\h_2$ de $\g$. On peut évidemment généraliser la construction en considérant des caractères. La quantification de Cattaneo-Felder définit donc deux algèbres de cohomologie $H^{\bullet}(\mu_0^{(1)}, \mathcal{A}_1)$ et $H^{\bullet}(\mu_0^{(2)}, \mathcal{A}_2)$. Dans \cite{CF1, CF2} Cattaneo-Felder  définissent une structure de bi-module sur une troisième algèbre de cohomologie.\\

La construction procède de la façon suivante.\\

 On fixe une décomposition de $\g$ compatible avec $\h_1$ et $\h_2$, c'est à dire que les variables porteront $4$ couleurs notées ici $(\pm, \pm )$. Le premier signe (resp. le second) vaut $+$ si la variable est tangente  et vaut $-$ si la variable est normale à $\h_1^\perp$ (resp. à $\h_2^\perp$).

On définit alors la fonction d'angle à $4$ couleurs dans le premier  quadrant  $0\leq \arg(z) \leq  \frac\pi{2}$ par la formule
\begin{equation}\phi_{\epsilon_1, \epsilon_2}(p,q)=\arg (p-q)+\epsilon_1 \arg
(p-\overline{q}) + \epsilon_2 \arg (p+ \overline{q})
+\epsilon_1\epsilon_2 \arg(p+q).\end{equation}
La fonction d'angle vérifie la propriété suivante :

-- lorsque $p ,\,q$ se concentrent sur l'axe horizontal
 les fonctions d'angles $d\phi_{\epsilon_1, \epsilon_2}(p,q)$
  tendent vers la $1$-forme d'angle

  \begin{equation}\mathrm{d}\phi_{\epsilon_1}(p,q)= \mathrm{d}\arg (p-q)+\epsilon_1 \mathrm{d}\arg
(p-\overline{q}), \end{equation}

-- lorsque $p ,\,q$ se concentrent sur l'axe vertical
 les fonctions d'angles $\mathrm{d}\phi_{\epsilon_1, \epsilon_2}(p,q)$
  tendent vers la $1$-forme d'angle

  \begin{equation}\mathrm{d}\phi_{\epsilon_2}(p,q)= \mathrm{d}\arg (p-q)+\epsilon_2 \mathrm{d}\arg
(p+\overline{q}).  \end{equation}

On dessine dans le premier  quadrant tous les diagrammes $\Gamma$ de Kontsevich colorés par les $4$ couleurs ci-dessus en plaçant les  sommets de première espèce dans le quadrant strict  et les sommets de deuxième espèce sur les axes.\\

En considérant les compactifications  des  configurations de points du premier quadrant modulo l'action du groupe des dilatations, on définit alors des variétés à coins compactes sur lesquelles on pourra intégrer les formes $\Omega_\G$ associées\footnote{A chaque arêtes  colorés est associée la différentielle d'une des $4$ fonctions d'angle ci-dessus.}.\\

Chaque graphe coloré $\G$ va définir un opérateur polydifférentiel, une fois que l'on aura placé  aux sommets de deuxième espèce (placés sur les axes) des fonctions\footnote{J'entends des fonctions avec composantes impaires.}. Le résultat est alors restreint à $(\h_1+\h_2)^\perp$. Les arêtes colorées $(-,\pm)$ qui arrivent sur l'axe horizontal définissent des formes d'angles triviales, on peut donc placer sur  \textbf{l'axe horizontal} une fonction de $\mathcal{A}_1$. De manière analogue on placera des fonctions de $\mathcal{A}_2$ sur \textbf{l'axe vertical}.  Enfin on place à l'origine une fonction de \[\mathrm{Poly}\Big((\h_1+\h_2)^\perp\Big)\otimes \bigwedge(\h_1^*\cap\h_2^*).\]
Cette algèbre est munie d'une différentielle $\mu^{(2,1)}$ correspondant aux contributions de tous les graphes colorés avec une arête sortant à l'infini  colorée par $\h_1^*\cap\h_2^*$, et  un seul sommet de deuxième espèce placé à l'origine.\\

En utilisant la formule de Stokes et en faisant l'inventaire des toutes les contributions, on montre le théorème  de compatibilité suivant
\begin{theo}[\cite{CF1, CT}] L'espace $H^{\bullet}(\mu_0^{(2)}, \mathcal{A}_2)$ agit par la gauche sur \break  $H^{\bullet}(\mu_0^{(2,1)}, \mathrm{Poly}\Big((\h_1+\h_2)^\perp\Big)\otimes \bigwedge(\h_1^*\cap\h_2^*) )$. L'espace $H^{\bullet}(\mu_0^{(1)}, \mathcal{A}_1)$ agit par la droite. On note $ \underset{1}\star$ l'action à droite (axe horizontal)  et $ \underset{2}\star $ l'action à gauche (axe vertical).
\end{theo}

\section{Applications en théorie de Lie}

On décrit maintenant les applications en théorie de Lie de la bi-quantification et du théorème de compatibilité.

\subsection{Description de l'algèbre de réduction}\label{sectionreduction}

On fixe un supplémentaire $\q$ de $\h$. On a donc une décomposition de l'algèbre enveloppante $$U(\g)=\beta(S(\q)) \oplus U(\g)\cdot \h_{-\lambda}$$ où $\beta$ désigne la symétrisation. On notera $\beta_\q$ l'application déduite de $S(\q)$ dans $U(\g)/U(\g)\cdot \h_{-\lambda}$.\\

On considère la bi-quantification
Cattaneo-Felder  pour le couple de
variétés co-isotropes $\h_{-\lambda}^\perp$ mis en position horizontale et
$0^\perp=\g^*$ mis en position verticale. L'espace de réduction associé à l'origine est tout simplement $S(\g)/S(\g)\cdot \h_{-\lambda}=\mathrm{Poly}(\h_{-\lambda}^\perp)$.  Comme dans \cite{Kont} on
considère
 pour le bi-vecteur  de Poisson, la moitié du crochet de Lie.\\

 \begin{figure}[h!]
\begin{center}
\includegraphics[width=12cm]{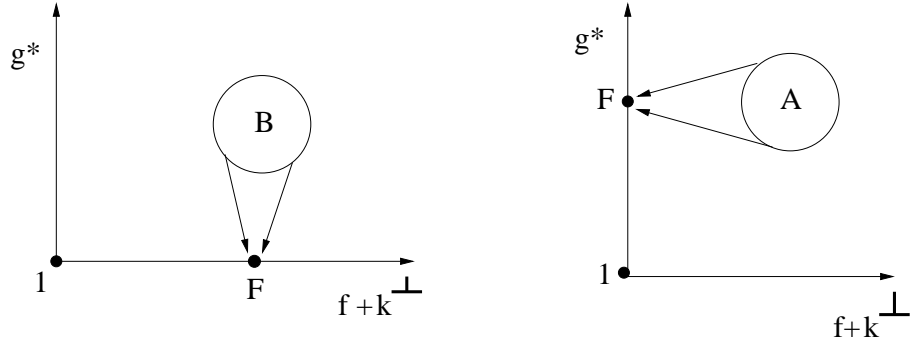}
\caption{\footnotesize Contributions des roues pures sur les axes}
\end{center}
\end{figure}

Pour $F\in S(\g)$ dans l'espace de réduction vertical,   $ F \underset{2}\star 1 \in S(\g)/S(\g)\cdot \h_{-\lambda}$  et $F\mapsto F \underset{2}\star 1 $ est un opérateur donné par diagrammes de Kontsevich \`a
$4$-couleurs. Cet opérateur est compliqué et n'est pas à coefficients
constants car les diagrammes avec ar\^etes
doubles colorées par $+\, + $ et $+ \, -$ ne sont pas nuls \textit{a priori}.

Lorsque $F \in S(\q)$ alors l'opérateur $A(F)=F\underset{2}\star 1$ est plus simple et correspond à l'exponentielle des contributions des roues pures sur l'axe vertical. On notera de même  $B(F)=1 \underset{1}\star F$ l'exponentielle des contributions des roues sur l'axe horizontal. On note $A(X)$ et $B(X)$ les symboles associés, c'est à dire pour $X\in \q$, $A_\q(X)=(e^X\underset{2}\star 1) e^{-X}$ et $B_\q(X)= (1 \underset{1}\star e^X)  e^{-X}$.\\

Je note pour $X\in \g$, $j(X)= \det_{\g}\Big(\frac{\sinh \frac{\ad X}2}{\frac{\ad
X}2}\Big)$ et on définit $J_\q(X)$ par la formule

\begin{equation}A_\q(X) J_\q^{1/2}(X)=B_\q(X) j^{1/2}(X).\end{equation}

\begin{theo}
L'application $\beta_\q \circ J_\q^{1/2}(\partial)$ définit un isomorphisme d'algèbres de l'algèbre de réduction $H^0(\mu_0)$  muni du produit $\underset{CF}\star$  sur l'algèbre $\big(U(\g)/U(\g)\cdot \h_{-\lambda} \big)^\h$.
\end{theo}
Ce théorème est démontré dans le cas des paires symétriques dans \cite{CT} et est étendu au cas général dans la thèse de mon étudiant P. Batakidis \cite{BP}. La formule de Stokes est encore l'argument essentiel. Lorsque $\q$ est invariant par action adjointe de $\h$ alors $J_\q$ est une fonction $\ad(\h)$-invariante. Dans le cas des paires symétriques on montre \cite{CT} \S 4.1 que $J_\q$ vaut $J(X)=\det_{\q}\Big(\frac{\sinh(\ad X)}{\ad X}\Big)$. Dans le cas des paires symétriques on retrouve le produit $\underset{Rou} \sharp$ de Rouvière \cite{Rou86, Rou90, Rou91, Rou94}.\\

On dispose donc d'une description complète  via la symétrisation des éléments de $\big(U(\g)/U(\g)\cdot \h_{-\lambda} \big)^\h$.

\subsection{Cas des paires symétriques}\label{pairesym}
Supposons dans cette section que $(\g, \sigma)$ est une paire symétrique, c'est à dire que $\sigma$ est une involution de Lie. On a donc une décomposition $\g=\k\oplus \p$ avec $[\k, \p]\subset  \p$ et $[\p, \p]\subset \k$. On s'intéresse à l'algèbre $\Big(U(\g)/U(\g)\cdot \k\Big)^\k$.  Le produit $\mu_1$ prend alors une forme remarquable. En effet pour $X, Y \in \p$ on a

\begin{equation}\mu_1(e^X, e^Y)=E(X, Y) e^{X+Y},\end{equation}
où $E(X, Y)$ est l'exponentielle des contributions des graphes de types roues\footnote{C'est à dire des roues colorées, attachées à des graphes de type Lie.}. Le produit de $\underset{CF}\star=\underset{Rou}\sharp$ se résume sur $H^0(\mu_0)=S(\p)^\k$, les éléments $\k$-invariants, en un opérateur  bidifférentiel à coefficient constant de symbole~$E$.\\

On peut déduire de l'analyse de la fonction $E$ des propriétés non triviales de  $\Big(U(\g)/U(\g)\cdot \k\Big)^\k$. En voici quelques unes (cf. \cite{CT} \S 3) :

 \paragraph{Symétrie :} Pour $X, Y\in \p$, on a $E(X, Y)=E(Y, X)$. On déduit la commutativité de  l'algèbre $\Big(U(\g)/U(\g)\cdot \k\Big)^\k$.\\

\paragraph{Radical résoluble :} Si  $X$ est dans le radical résoluble\footnote{Le plus grand idéal résoluble de $\g$, il est $\sigma$-stable.} de $\g$, alors $E(X,Y)=1$ pour tout $Y\in \p$. On en déduit que $F\underset{CF}\star G= FG$ si $F \in S(\p\cap J)^\k$ avec $J$ un idéal résoluble.\\

\paragraph{Paires symétriques d'Alekseev-Meinrenken:} Supposons que $(\g, \sigma)$ soit munie d'une forme bilinéaire non dégénérée\footnote{On dit que $(\g, \sigma)$ est une paire symétrique quadratique.} et $\sigma$-anti-invariante, alors on a $E(X, Y)=1$.\\

\paragraph{Double quadratique:} Supposons que $\g$ soit une algèbre de Lie quadratique, c'est à dire munie d'une forme bilinéaire non dégénérée. Considérons la paire symétrique double $(\g_{double}, \sigma)=(\g\times\g, \sigma)$ avec $\sigma(X, Y)=(Y, X)$. Alors on a $E_{double}=1$. On conjecture que cette propriété est vraie même si $\g$ n'est pas quadratique, ce qui a d'intéressantes conséquences (voir la fin de l'article).

\subsection{Opérateurs différentiels en coordonnées expo\-nentielles}

On reprend dans cette section les notations de l'introduction. On peut raffiner le théorème précédent en donnant l'écriture en coordonnées exponentielles des opérateurs différentiel invariants du fibré $\mathcal{L}_\lambda$. On considère toujours le  diagramme de bi-quantification du \S \ref{sectionreduction}. Notons $\Exp$ l'application exponentielle de $\q$ sur $G/H$. On travaille au voisinage de $0\in \q$. Soit $\varphi \in\G(\mathcal{L}_\lambda)$, alors $\varphi$ est une fonction sur $G$ telle que $\varphi(gh)=\varphi(g)\lambda(h)$. Notons $\phi \in \mathcal{C}^\infty(\q)$ définie au voisinage de $0$ par $$\phi(X)=\frac{J_\q^{1/2}(X)}{B(X)} \times \varphi(\exp_G(X)).$$Le facteur représente une sorte de jacobien. Pour $u\in \mathcal{D}_\lambda$ on note $D_u^{\Exp}$ l'opérateur en coordonnées exponentielles défini par $$D_u^{\Exp}(\phi)(X)= \frac{J_\q^{1/2}(X)}{B(X)} \times D_u(\varphi)(\exp_G(X)).$$

Pour $X, Y \in \q$ on note $Q(X, Y)\in \q$ et $H(X, Y)\in \h$ les composantes exponentielles au voisinage de~$0$ :
\begin{equation}\exp_G(X)\exp_G(Y)= \exp_G(Q(X, Y))\exp_G(H(X, Y)).\end{equation}
On peut étendre au cas des sous-algèbres la formule de \cite{CT} \S 4.3. Pour $R$ dans l'algèbre de réduction $H^{0}(\mu_0)$ on a :

\begin{equation}e^X\underset{1}\star R=\frac{J_\q^{1/2}(X)}{B(X)} \times R(\partial_Y) \Big( J_\q^{1/2}(Y) \times \frac{B(Q(X,Y))}{J_\q^{1/2}(Q(X,Y))} e^{Q(X,Y)}e^{\lambda(H(X, Y))} \Big)\mid_{Y=0}.\end{equation}

\begin{theo}[\cite{CT}]\label{theoopd} Pour $u=\beta\Big(J_\q^{1/2}(\partial)(R)\Big)$, l'opérateur différentiel $D_u^{\Exp}$ en coordonnées exponentielles s'exprime par la formule $e^X\underset{1} \star R$.
\end{theo}
Remarquons que l'expression est valable au voisinage de $0$ en posant $\epsilon=1$. La convergence des coefficients se démontre comme dans \cite{ADS, AST}. Cette expression résout de manière satisfaisante un problème ancien de M. Duflo \cite{duflo-japon}.\\

\noindent \textbf{Remarque :} Ce théorème devrait avoir des conséquences intéressantes dans le cas des espaces symétriques hermitiens. En effet on peut utiliser la réalisation d'Harish-Chandra pour écrire les opérateurs différentiels invariants\footnote{Les formules proposées dans \cite{Nomu} ne sont pas correctes.} \cite{Shi, Nomu}.

\subsection{Construction de caractères pour l'algèbre $\big(U(\g)/U(\g)\cdot \h_{-\lambda} \big)^\h$}

On reprend dans cette section les notations de l'introduction. On cherche des caractères pour l'algèbre $\Big(U(\g)/U(\g)\cdot \h_{-\lambda}\Big)^\h$. Afin que ceci soit intéressant il faut que cette dernière soit commutative.\\

\subsubsection{L'hypothèse lagrangienne} Plaçons nous dans l'hypothèse lagrangienne que nous avons décrite dans l'introduction. On suppose que pour $f\in \h_{-\lambda}^\perp$ générique, l'espace  $\h\cdot f$ (action coadjointe) est lagrangien dans $\g\cdot f$. Rappelons que $\g\cdot f$ est toujours un espace symplectique muni de la forme de Kostant-Souriau. Il s'identifie à l'espace symplectique $(\g/\g(f), B_f)$ avec $\g(f)=\{X\in \g, \; X\cdot f=0\;  \}$ et $B_f(X, Y)=f[X, Y]$.\\

 Comme on l'a dit l'algèbre de Poisson $\mathcal{S}_\lambda= \Big(S(\g)/S(\g)\cdot \h_{-\lambda}\Big)^\h$ est alors  Poisson commutative, mais aussi l'algèbre des fractions invariantes\footnote{On notera qu'il faut prendre les invariants du corps des fractions et non pas le corps des fractions des invariants.} $\Big(\mathrm{Frac}(\h_{-\lambda}^\perp)\Big)^\h$. En fait la commutativité de cette dernière est équivalente à l'hypothèse lagrangienne.\\

 On suppose de plus que les $H$-orbites génériques sont polarisables. C'est à dire que pour $f\in \h_{-\lambda}^\perp$ générique, il existe une polarisation $\b$ en $f$ (sous-algèbre isotrope et de dimension maximale parmi les sous-espaces isotropes). On a $f[\b, \b]=0$ et $\b/\g(f)$ lagrangien dans $\g/\g(f)$. Cette hypothèse est au c\oe ur de la quantification géométrique.\\

  Dans ces conditions on montre facilement que $(\h\cap \b)\cdot f= (\h+\b)^\perp$. Cette hypothèse va être cruciale.\\
\subsubsection{Construction d'un caractère}
  Dans \cite{CT} \S 6, on construit le diagramme de bi-quantification en plaçant sur l'axe horizontal $\h_{-\lambda}^\perp$ et sur l'axe vertical $f + \b^\perp$. Par construction il faut donc choisir une décomposition compatible de $\g$. En particulier on fixe un supplémentaire $\q$ de $\h$ en position d'intersection normale avec $\b$ c'est à dire que l'on a $$\b=\b\cap \h \oplus \b\cap \q.$$
  D'après \S \ref{sectionexemple} on sait que l'algèbre de réduction verticale est réduite à $\mathbb{C}$ et l'hypothèse $(\h\cap \b)\cdot f= (\h+\b)^\perp$ implique facilement que  l'algèbre de réduction associée à l'origine est aussi $\mathbb{C}$. On en déduit une action à droite de $H^0(\mu_0) \sim \Big(U(\g)/U(\g)\cdot \h_{-\lambda}\Big)^\h$ sur $\mathbb{C}$. On a donc construit un caractère de cette algèbre.\\

\begin{figure}[h!]
\begin{center}
\includegraphics[width=5cm]{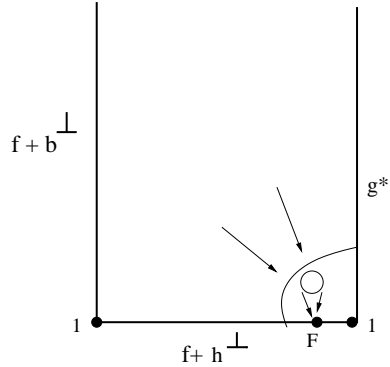}
\caption{\footnotesize Calcul du caractère}\label{caractere.eps}
\end{center}
\end{figure}

   Dans \cite{CT} \S 6 on développe une théorie des diagrammes à $8$-couleurs dans une bande qui nous permet d'interpoler deux situations de bi-quantification et nous permet de calculer le caractère (cf. Fig. \ref{caractere.eps}). On peut déplacer la position de $F$ le long du bord horizontal. Les positions limites aux coins donnent les informations recherchées. Notre méthode utilise la formule de Stokes.

  \begin{theo} Sous les hypothèses lagrangiennes ci-dessus et l'existence de polarisation, l'application  $$u \in  \Big(U(\g)/U(\g)\cdot \h_{-\lambda}\Big)^\h \mapsto B_\q(\partial)J_\q^{-1/2}(\partial)(\beta_\q ^{-1}(u)) (f)$$ est le caractère construit par le diagramme de bi-quantification.
 \end{theo}
 Ce théorème est démontré dans le cas des paires symétriques dans \cite{CT} \S 6 et étendu au cas général dans la thèse de mon étudiant P. Batakidis \cite{BP}. Remarquons que dans le cas des paires symétriques résolubles, on retrouve directement la formule de Rouvière, car $B_\p(X)=1$ et $J_\p(X)=\det_{\q}\Big(\frac{\sinh(\ad X)}{\ad X}\Big)$.\\

 On en déduit facilement en regardant le terme dominant le corollaire suivant

 \begin{cor} Sous les hypothèses lagrangiennes ci-dessus, l'algèbre \break $\big(U(\g)/U(\g)\cdot \h_{-\lambda} \big)^\h$ est commutative.
 \end{cor}
En particulier dans le cas nilpotent, les orbites sont  toujours polarisables\footnote{On peut prendre une polarisation construite par M. Vergne.}, l'hypothèse lagrangienne est équivalente au fait que la multiplicité de la représentation $L^2(G/H, \chi)$ est finie\footnote{Ici $\chi$ est le caractère de différentielle $i\lambda$.}. On retrouve alors un théorème de Corwin \& Greenleaf \cite{CG}.  On en déduit dans le cas nilpotent,  que si l'algèbre  $S_\lambda$ est Poisson commutative alors  $\big(U(\g)/U(\g)\cdot \h_{-\lambda} \big)^\h$ est commutative. Réciproquement d'après le lemme \ref{com},  si $\big(U(\g)/U(\g)\cdot \h_{-\lambda} \big)^\h$ est commutative alors la limite classique $LC(H^0(\mu_0))$ est Poisson commutative. Bien évidemment si cette dernière algèbre de Poisson est assez grosse\footnote{Une conjecture raisonnable semble que le corps des fractions de la limite classique contient $\mathcal{S}_\lambda$.} on pourrait en déduire alors que $S_\lambda$ est aussi Poisson commutative.

\subsubsection{Comparaison avec le vecteur de Penney}

Sous les hypothèses lagrangiennes ci-dessus, la méthode des orbites fournit dans les bons cas un vecteur distribution semi-invariant (dit vecteur de Penney \cite{To3, To4, FLMM, Li1, Li2, Li3}). Dans les bons cas  ce vecteur est aussi un vecteur propre pour $\mathcal{D}_\lambda$, ce qui fournit aussi un caractère de cette algèbre.\\

Décrivons le vecteur de Penney dans la situation unimodulaire et $\lambda =0$ pour simplifier l'exposé. On note $B$ un groupe de Lie connexe d'algèbre de Lie $\b$. On suppose que $\chi_f(\exp_B(X))=e^{if(X)}$ définit bien un caractère\footnote{Ces hypothèses d'intégrabilité compliquent la théorie en général.} de $B$. On note $\d_{G, H}$ et $ \d_{B, B\cap H }$ des mesures invariantes \footnote{En général de telle mesure n'existe pas, il faut alors travailler sur des fibrés en droite, cf. Introduction.} sur $G/H$ et $B/B\cap H$.

Le vecteur de Penney est défini par la fonction généralisée

$$\Phi d_{G,H}\longmapsto j_*(\Phi d_{G,H})=\int_{B/B\cap H}
\Phi(b)\chi_f(b)^{-1} d_{B, B\cap K}(b).
$$

Sous l'hypothèse lagrangienne c'est
une section g\'en\'eralis\'ee propre \cite{To3} sous l'action des op\'erateurs
diff\'erentiels invariants $D_u$ pour $u\in (U(\g)/U(\g)\cdot
\h)^{\h}$. Ceci fournit donc un
caract\`ere de cette alg\`ebre not\'e $u \longrightarrow\lambda_{f,\b}(u)$.\\

En utilisant  le théorème~\ref{theoopd} on montre que ce caractère vaut aussi la transformée de Fourier du  symbole transverse et coïncide (dans les bons cas) alors avec celui construit par le diagramme de bi-quantification\footnote{Le cas nilpotent sera traité dans la thèse de P. Batakidis \cite{BP}.}.

\subsection{Dépendance par rapport au supplémentaire et applications }
On étudie maintenant la dépendance du produit $\underset{CF}\star$ par rapport au choix du supplémentaire $\q$.

\subsubsection{Formule de changement de base}

On fixe un supplémentaire $\q_0$
de $\h$. On choisit  $(e_i)$ une base de $\g$ adaptée à la décomposition $\h\oplus \q_0$, c'est à dire une base $(K_i)_i$ de $\h$ et une base $(P_a)_a$ de~$\q_o$. Faisons choix d'un autre supplémentaire $\q_1$ dont on fixe une base
$(Q_a)_a$. Sans perte de généralité on peut supposer que la matrice de passage est de la forme
\begin{equation}\nonumber
\mathbb{M}=\left(\begin{array}{cc}
 I & \mathbb{D}\\
  0 & I
\end{array}\right)
\end{equation}

Notons $\mathbb{D}=[V_1, \ldots, V_ p]$ les
colonnes de la matrice $\mathbb{D}$ et $V_i\in \h$.

On écrit le bi-vecteur $\pi$ dans les deux décompositions et on applique la procédure de transformée de Fourier impaire dans les directions normales. On trouve alors deux poly-vecteurs $\widehat{\pi}$ et $\widehat{\pi}^{(1)}$ sur la  variété intrinsèque  $\h^\perp \oplus \mathrm{\Pi}\h$. La relation entre les deux  poly-vecteurs est la suivante (cf. \cite{CT} \S 1.5).
Considérons
\[\pi_\mathbb{M}=\mathbb{M}^{-1}[\mathbb{M}e_i, \mathbb{M}e_j]\partial_{e_i^*}\wedge
\partial_{e_j^*}.\]
Alors on aura \[\widehat{\pi}^{(1)}=\widehat{\pi_\mathbb{M}}\]où le membre de droite est la transformée de Fourier partielle  impaire pour la première décomposition. Par ailleurs le champ de vecteurs sur $\g^*$ défini par  $v= -V_a
\partial_{P_a^*}$ vérifie\footnote{On note $[\bullet, \bullet]_{S}$ le crochet de
Schouten-Nijenhuis.}  $[v,v]_{S}=0$ et on a la relation

\begin{equation}\pi_\mathbb{M}= e^{\ad v} \cdot \pi= \pi + [v, \pi]_{S}+\frac {1}{2}[v,[v,\pi]]_{S}.\end{equation}
L'action du champ $-v$ sur le bivecteur $\pi$ correspond au
changement de supplémentaire.

\subsubsection{Contrôle de la déformation: l'élément de jauge}
Soit $t$ un paramètre réel. Je note \[\pi_t= e^{t\,  \ad v}\cdot \pi =\pi + t[v, \pi] + \frac {t^2}2 [v,[v,\pi]]]\] et $\widehat{\pi_t}$ sa transformée de Fourier partielle pour la décomposition $\g=\h\oplus \q_0$.\\

Pour $t=0$ on trouve $\pi$ et pour $t=1$ on trouve $\pi_\mathbb{M}$.  On applique le $L_\infty$ quasi-isomorphisme du théorème \ref{theoKont} à la super-variété $\h^\perp\oplus\mathrm{\Pi} \h$. Alors  $\mu_t$ définie par
\[ \mu_t=\mathcal{U}(
e^{\widehat{\,\pi_t}})=m+\sum_{n\geq 1} \frac {\epsilon^n}{n!}U_n
\left(\widehat{\,\pi_t}, \ldots, \widehat{\,\pi_t}\right)\]est une structure $A_\infty$.\\

 La dérivée $\d\mathcal{U}_{\widehat{\,\pi_t}}$ au point
$\widehat{\,\pi_t}$ est un morphisme de complexes. On a donc en
dérivant, une équation différentielle linéaire  :

\begin{equation}\frac{\partial \mu_t}{\partial t}= \epsilon
\,\d\mathcal{U}_{\widehat{\,\pi_t}} \left({[\widehat{\,v},
\widehat{\,\pi_t}]}_{SN}\right)={[\d\mathcal{U}_{\widehat{\,\pi_t}}
\left(\widehat{\,v}\right), \mu_t]}_{G} \; .\end{equation}
On peut alors traduire cette formule en terme de diagrammes de Kontsevich colorés comme dans \S \ref{sectionrappels}. On placera aux sommets de première espèce le bi-vecteur $\pi_t$ et une fois le vecteur $v$. Les  sommets terrestres reçoivent des fonctions de $\Big(S(\g)/S(\g)\cdot\h\Big)\otimes
\bigwedge \h^*$.\\

\begin{figure}[h!]
\begin{center}
\includegraphics[width=6cm]{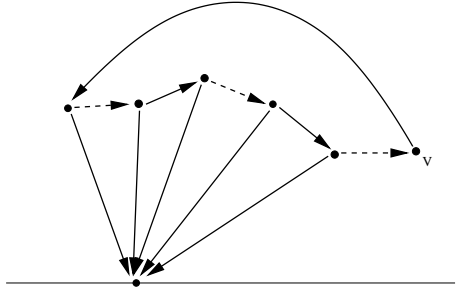}
\caption{\footnotesize Bernoulli fermé par $v$}\label{HC1.eps}
\end{center}
\end{figure}

\begin{figure}[h!]
\begin{center}
\includegraphics[width=6cm]{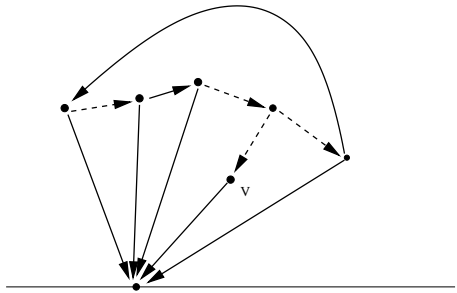}
\caption{\footnotesize Roue pure attachée  à $v$}\label{HC2.eps}
\end{center}
\end{figure}

\begin{figure}[h!]
\begin{center}
\includegraphics[width=8cm]{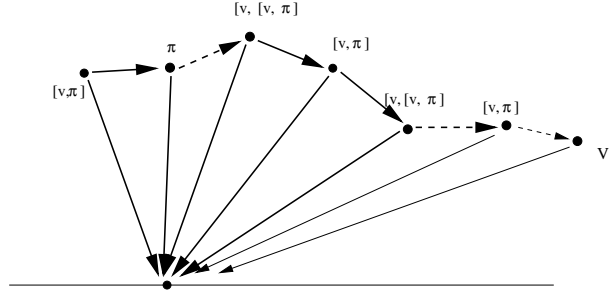}
\caption{\footnotesize Graphe de type Bernoulli non
fermé}\label{bernoulliv.eps}
\end{center}
\end{figure}

La  différentielle $(\mu_t)_0$, composante de degré~$0$ de $\mu_t$, vérifie en particulier l'équation différentielle

\begin{equation}\label{jauge}\frac{\partial (\mu_t)_0}{\partial t}= [ \left(DU_{\widehat{\pi_t}}(\widehat{v})\right)_0,(\mu_t)_0 ].\end{equation}
Cette formule dit que les différentielles  $(\mu_t)_0$ sont conjuguées par un élément de type groupe\footnote{C'est la résolvante de l'équation différentielle.}: c'est l'élément de jauge. Pour le décrire, il suffit d'analyser tous les graphes qui interviennent dans $\left(DU_{\widehat{\pi_t}}(\widehat{v})\right)_0$. C'est ce que l'on fait dans \cite{CT}~\S~5.5. Ce sont les graphes des figures Fig.~\ref{bernoulli2.eps}, Fig.~\ref{rouebernoulli.eps} et Fig.~\ref{roueinfini.eps} où l'arête $\infty$ dérive le sommet attaché au vecteur $v$. L'arête issue de $v$ va soit sur la racine du graphe, soit sur le sommet terrestre. Il y a donc $4$ types de graphes donnés par les figures suivantes (Fig. \ref{HC1.eps}, Fig. \ref{HC2.eps}, Fig. \ref{HC2.eps}  et Fig. \ref{bernoulliv.eps}).\\

\begin{figure}[h!]
\begin{center}
\includegraphics[width=8cm]{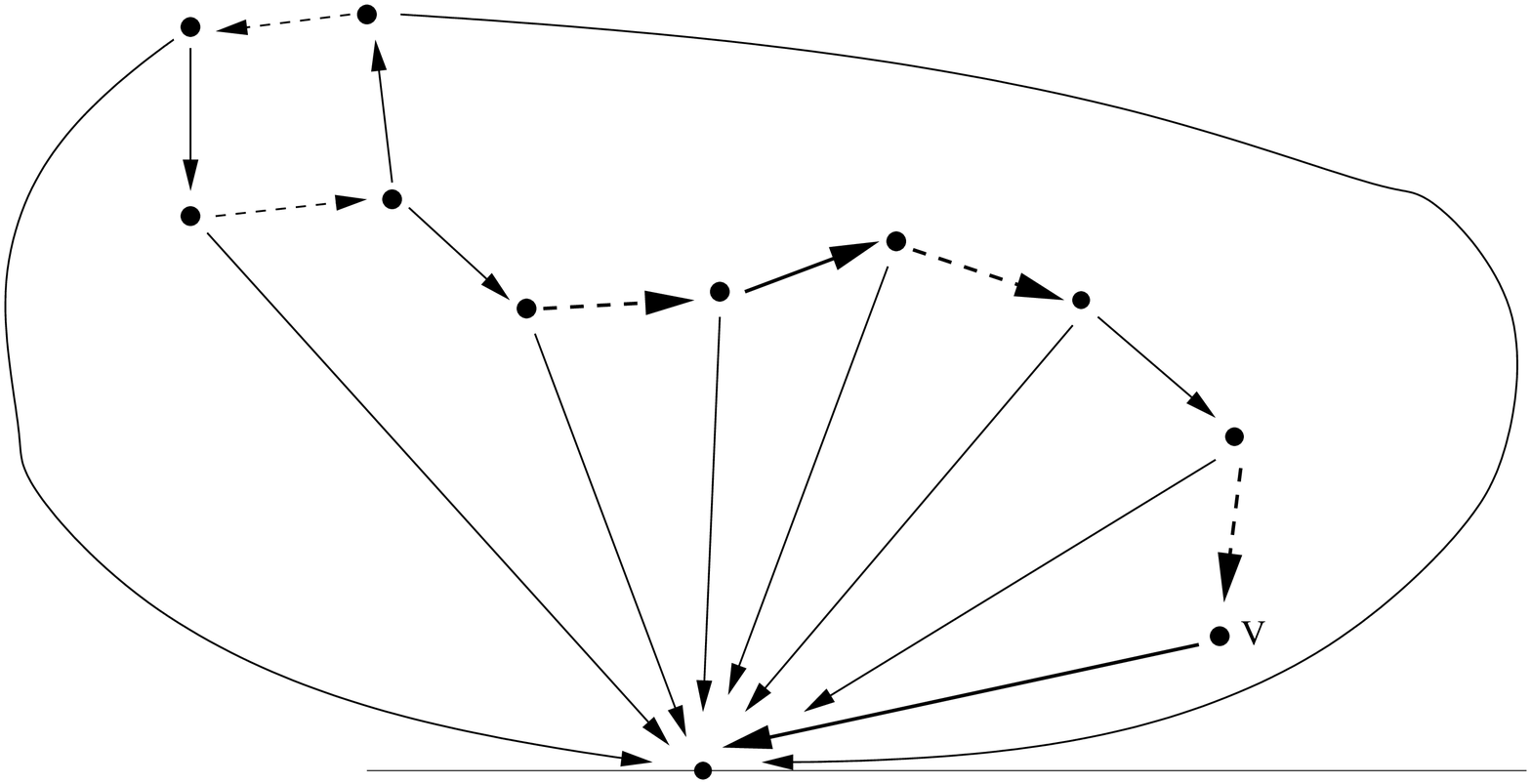}
\caption{\footnotesize Roue attachée à un Bernoulli attaché à $v$}\label{HC3.eps}
\end{center}
\end{figure}

De même,  la composante de degré un $(\mu_t)_1$ vérifie l'équation

\begin{equation}\label{mu1}\frac{\partial (\mu_t)_1}{\partial t}= [ \left(DU_{\widehat{\pi_t}}(\widehat{v})\right)_0,(\mu_t)_1 ]+[ \left(DU_{\widehat{\pi_t}}(\widehat{v})\right)_1,(\mu_t)_0 ] .\end{equation}

Cette équation exprime alors que l'élément de jauge définit en cohomologie un isomorphisme d'algèbres. En particulier, résoudre l'équation (\ref{jauge}) permet de trouver explicitement l'entrelacement des star-produits pour deux choix de supplémentaires.

\subsubsection{Exemple du double}
Pour illustrer les conséquences de nos théories, examinons le cas des paires symétriques $(\g\times \g, \sigma)$ avec $\sigma(X, Y)=(Y,X)$. On a alors  $\k=\{(X, X), \; X\in \g \}$. Ce sont les doubles.

\paragraph{Déformation du supplémentaire :}  On dispose de $3$ supplémentaires invariants naturels, $$\g_{-}=\{(0, 2X)\;, X\in \g \} \quad \p=\{(X, -X)\;, X\in \g \} \quad \g_{+}=\{(2X, 0)\;, X\in \g \}.$$
L'interpolation est donnée par la famille à paramètre de sous-espaces invariants   $$\{((1+t)X,(t-1)X)\;, X\in \g \}.$$Le champ de vecteur $v= -V_a
\partial_{P_a^*}$ est donc l'application linéaire qui transforme $(X, -X)\in \p \mapsto -(X, X) \in \k$.  Si on note $(e_i)$ une base de $\g$, $K_i=(e_i, e_i)$ et $P_j=(e_j, -e_j)$ alors le bi-vecteur $\pi_t$ prend une forme assez simple

\begin{multline}\pi_t=2[K_i, P_j] \partial_{K_i^*}\wedge\partial_{P_j^*} + [K_i, K_j] \partial_{K_i^*}\wedge\partial_{K_j^*} + \\\left([K_i, K_j](1-t^2) + 2t [K_i, P_j]\right) \partial_{P_i^*}\wedge\partial_{P_j^*}
\end{multline}

\paragraph{Élément de jauge :} Comme  $(X, -X)$ et $(X, X)$ commutent l'opérateur $\left(DU_{\widehat{\pi_t}}(\widehat{v})\right)_0$ va se simplifier.
En effet les  opérateurs différentiels correspondant aux graphes des Fig.~\ref{HC3.eps} et Fig.~\ref{bernoulliv.eps} seront nuls. Il ne reste alors que les graphes Fig.~\ref{HC1.eps} et Fig.~\ref{HC2.eps}. Le coefficient vérifie une condition de symétrie, il est nul si le nombre de sommets attachés à $\pi_t$ est impair.
Au final ces diagrammes  correspondent à des opérateurs différentiels de symboles  $$P_{2n}(t) \tr_\p \big(\ad (X, -X)^{2n}\big)= P_{2n}(t) \tr_\g (\ad X)^{2n} $$ avec $P_{2n}(t)$ un polynôme  en  $t$ de degré inférieur à $2n$.
Ce polynôme dépend de la couleur dans la roue et des coefficients de Kontsevich associés, il s'écrit sous la forme $$P_{2n}(t)=\sum w_\G P_\G(t)$$
où $\G$ décrit tous les graphes Fig.~\ref{HC1.eps} et Fig.~\ref{HC2.eps} avec $2n$ sommets associés $\pi_t$. Ces polynômes sont sans doute liés aux polynômes de Bernoulli, mais nous n'avons pas pu le vérifier.\\

L'opérateur $\left(DU_{\widehat{\pi_t}}(\widehat{v})\right)_0$  commute à l'action adjointe, donc  l'équation de  la différentielle (\ref{jauge}) se résout simplement, on trouve $(\mu_t)_0=\epsilon \d_{CE}$, ce que l'on savait déjà par ailleurs.\\

Enfin les opérateurs $\left(DU_{\widehat{\pi_t}}(\widehat{v})\right)_0$ forment une famille commutative en $t$ donc l'élément  de jauge est donné par la résolvante $$\phi_t := \exp\Big(\int_0^t \left(DU_{\widehat{\pi_s}}(\widehat{v})\right)_0 \d s.\Big)$$
C'est un opérateur universel de symbole $$\phi_t(X):=\exp\left(\sum\limits_{n>0} Q_{2n+1}(t) \tr_\g (\ad X)^{2n}\right) .$$\\

\paragraph{Interpolation du produit :} Regardons le terme de degré~$1$ pour la structure $A_\infty$. L'espace $\h^\perp$ s'identifie aux couples $(f, -f)$ avec $f\in \g^*$. Les formes linéaires $$(2X, 0)\quad (X, -X) \quad (0,2X)$$sont identiques sur $\h^\perp$. On la note $\widehat{X}$.\\

Pour $t=0$, c'est la situation des paires symétriques on a donc d'après~\S~\ref{pairesym}
$$\mu_{t=0}\left(e^{\widehat{X}}, e^{\widehat{Y}}\right)=E_{double}((X,-X) , (Y,-Y) )e^{\widehat{X+Y}}.$$

Pour $t=1$, c'est la situation des algèbres de Lie, on a donc \cite{AST}
$$\mu_{t=1}\left(e^{\widehat{X}}, e^{\widehat{Y}}\right)=D(2X,2Y)e^{\widehat{\frac12 Z(2X, 2Y)}},$$ avec $Z(X, Y)$ la formule de Campbell-Hausdorff et $$D(X,Y)=\frac{j_\g^{1/2}(X)j_\g^{1/2}(Y)}{j_\g^{1/2}(Z)}$$la fonction de densité de Duflo. Rappelons que $j_\g(X)= \det_{\g}\Big(\frac{\sinh \frac{\ad X}2}{\frac{\ad
X}2}\Big)$.

Pour $t=1$, il est facile de calculer la fonction $J_{\g_+}$ que nous avons introduite en \S~\ref{sectionreduction}. En effet, les roues $A, B$ n'ont qu'une seule couleur $(+,+)$. Par symétrie on trouve $A=B$.  On a donc $$J_{\g_+}((2X,0))= j_\g(2X)=\det_{\g}\Big(\frac{\sinh \ad X}{\ad X}\Big)=J_\p(X, -X).$$

La fonction $J$ est donc la même pour ces deux choix de supplémentaires.

\begin{prop} L'élément de jauge $\phi_1$ vaut $1$. Dans le cas du double, le star-produit $\underset{CF}\star$  est trivial sur $S(\p)^k$, les éléments $\k$-invariants de $S(\p)$.  \end{prop}
\noindent \textit{Preuve:} Cet élément est un isomorphisme pour les algèbres de réduction. Pour $t=0$ on trouve $S(\p)^\k$ muni du produit de Rouvière qui vaut $\underset{CF}\star$. Pour $t=1$ c'est le produit de Duflo-Kontsevich, c'est à dire la multiplication standard sur les invariants.  Or pour le double quadratique la fonction $E_{double}$  vaut $1$ \cite{CT}. Donc $\phi_1$ est un isomorphisme d'algèbres de $S(\g)^\g$ pour toute algèbre quadratique. C'est donc l'identité car $\tr_\g (\ad X)^{2n}$ n'agit pas comme une dérivation universelle. Comme $\phi_1$ est universel, c'est toujours $1$.  On en déduit que l'action de $E_{double}$ sur $S(\p)^k$ est triviale, même si $\g$ n'est pas quadratique. \fin\\

\paragraph{Conjecture de Kashiwara-Vergne :}
Plus généralement, comme les formules de quantification dans le cas linéaire sont toujours des exponentielles on aura pour la composante de degré~1 de notre structure $A_\infty$ : $$\mu_{t}\left(e^{\widehat{X}}, e^{\widehat{Y}}\right)=E_{t}(X,Y)e^{\widehat{Z_t(X,Y)}}$$avec $Z_t(X,Y)$ une série de Lie formelle en $X, Y$ à coefficients polynomiaux en~$t$.

On dispose donc d'une déformation de la formule de Campbell-Hausdorff en la loi additive. On peut alors traduire l'équation (\ref{mu1}) sur la déformation. On compense déjà le terme $[ \left(DU_{\widehat{\pi_t}}(\widehat{v})\right)_0,(\mu_t)_1 ]$ en conjuguant par $\phi_t(X)$. Posons donc
$$\widehat{E_t(X, Y)}= E_t(X, Y) \frac{\phi_t(X)^{-1}\phi_t(Y)^{-1}}{\phi_t(Z_t)^{-1}}.$$

En examinant ce qu'est le terme $\left(DU_{\widehat{\pi_t}}(\widehat{v})\right)_1$ on se convainc sans difficulté, comme dans \cite{To1} que l'on est en train de calculer une différentielle en $X$ et en $Y$ de la fonction $Z_t(X, Y)$. On trouve   un contrôle à la Kashiwara-Vergne \cite{KV} de la déformation $Z_t(X, Y)$ (voir \cite{To5} pour un résumé de les méthodes de Kashiwara-Vergne). On a donc montré la théorème  suivant.

\begin{theo} La déformation du supplémentaire produit une déformation de Kashiwara-Vergne, c'est à dire qu'il existe des séries de Lie sans termes constants $(F_t(X, Y), G_t(X, Y))$ à coefficients polynomiaux en $t$,  telles que

\begin{equation}\partial_t Z_t(X, Y) = [X, F_t(X, Y)]\cdot \partial_X Z_t(X, Y) + [Y, G_t(X, Y)]\cdot \partial_Y Z_t(X, Y)\end{equation}
\begin{multline}\partial_t \widehat{E_t(X, Y)}= \Big([X, F_t(X, Y)]\cdot \partial_X  + [Y, G_t(X, Y)]\cdot \partial_Y \Big) \widehat{E_t(X, Y)}  \; +  \\ \widehat{E_t(X, Y)}\tr_\g \big(\partial_X F_t \circ \ad X + \partial_Y G_t\circ \ad Y\big). \end{multline}
\end{theo}

\begin{cor} Si la conjecture $E_{double}=1$ est vraie,  alors la déformation du supplémentaire démontre la conjecture de Kashiwara-Vergne.
\end{cor}
\noindent \textit{Preuve: } En effet, si $E_{double}=1$ alors le théorème précédent fournit une déformation à la Kashiwara-Vergne qui à les bonnes conditions limites. C'est à dire pour $t=0$ on a le produit $X+Y$ est pour $t=1$ le produit de Duflo. Comme dans \cite{AM2, AT, To5} on construit alors une solution de Kashiwara-Vergne.

\fin \\

\noindent \textbf{Remarque finale :} Dans \cite{Kont, To1} les arguments d'homotopie se fondent sur la déformation géométrique réelle des coefficients $w_\G$. En regardant les algèbres de Lie comme des paires symétriques, on construit ici une déformation polynomiale des coefficients, ce qui est bien meilleur. On peut donc espérer que notre déformation est rationnelle.

\end{document}